\documentclass{elsart}
\usepackage[latin1]{inputenc}
\RequirePackage{amsfonts,amsmath,amssymb}
\usepackage{dsfont}

\newtheorem{theo}{Theorem}

\newtheorem{propo}{Proposition}[section]

\newtheorem{coro}{Corollary}[theo]

\numberwithin{equation}{section}
\newcommand{\dd}{\delta}
\renewcommand{\e}{\epsilon}
\newcommand{\alp}{\alpha}
\newcommand{\lab}{\lambda}

\newcommand{\sig}{\sigma}
\newcommand{\tht}{\theta}

\newcommand{\Gam}{\Gamma}
\newcommand{\Sig}{\Sigma}

\newcommand{\YYN}{(Y_1,\ldots,Y_n)}

\newcommand{\UUN}{(U_1,\ldots,U_n)}

\newcommand{\indep}{{\bot}\kern-0.9em{\bot}}
\newcommand{\wap}{(\Omega,\mathcal{A},\rm I\kern-2pt P)}
\newcommand{\Cov}{\mathrm{Cov}}
\newcommand{\Var}{\mathrm{Var}}
\newcommand{\sli}{\sum\limits}
\newcommand{\sliin}{\sum\limits_{i=1}^n}
\newcommand{\sliim}{\sum\limits_{i=1}^m}

\newcommand{\ili}{\int}
\newcommand{\proli}{\prod\limits}
\newcommand{\bculi}{\bigcup\limits}
\newcommand{\bcali}{\bigcap\limits}

\newcommand{\limk}{\lim_{k\rightarrow\infty}\;}
\newcommand{\lsn}{\limsup_{n\rightarrow\infty}}
\newcommand{\lin}{\liminf_{n\rightarrow\infty}}

\newcommand{\rar}{\rightarrow}

\newcommand{\cvps}{\rightarrow_{a.s.}\;}

\newcommand{\kif}{k\rightarrow\infty}
\newcommand{\nif}{n\rightarrow\infty}

\newcommand{\aoo}{\Big\{}
\newcommand{\aff}{\Big\}}
\newcommand{\coo}{\Big [}
\newcommand{\cff}{\Big]}
\newcommand{\poo}{\Big (}
\newcommand{\pff}{\Big)}
\newcommand{\po}{\big (}
\renewcommand{\pf}{\big)}
\newcommand{\ao}{\big \{}
\newcommand{\af}{\big \}}

\newcommand{\pooo}{\bigg (}
\newcommand{\pfff}{\bigg)}
\newcommand{\aooo}{\bigg \{}
\newcommand{\afff}{\bigg \}}

\newcommand{\poooo}{\Bigg (}
\newcommand{\pffff}{\Bigg)}

\newcommand{\coooo}{\Bigg [}
\newcommand{\cffff}{\Bigg ]}
\newcommand{\mmi}{\mid\mid}

\newcommand{\Mid}{\Big |}
\newcommand{\Mmi}{\Big|\Big|}

\newcommand{\mmid}{\mid\mid_{[0,1]^d}}
\newcommand{\Mmid}{\Big | \Big |_{[0,1]^d}}

\newcommand{\nk}{n_k}
\newcommand{\nkm}{n_{k-1}}

\newcommand{\ank}{a_{n_k}}
\newcommand{\ankm}{a_{n_{k-1}}}

\newcommand{\EEE}{\mathbb{E}}
\newcommand{\NNN}{\mathbb{N}}
\newcommand{\PPP}{ \mathbb{P}}
\newcommand{\CCC}{\mathbb{C}}

\newcommand{\RRR}{\mathbb{R}}

\newcommand{\QQQ}{\mathbb{Q}}

\newcommand{\AAA}{\mathcal{A}}

\newcommand{\FF}{\mathcal{F}}
\newcommand{\TT}{\mathcal{T}}
\newcommand{\GG}{\mathcal{G}}

\newcommand{\SSS}{\mathcal{S}}
\newcommand{\PP}{\mathcal{P}}

\newcommand{\NN}{\mathcal{N}}

\newcommand{\DD}{\mathcal{D}}

\newcommand{\wt}{\widetilde}

\newcommand{\wtC}{\widetilde{C}}
\newcommand{\wtD}{\widetilde{D}}

\newcommand{\lb}{\newline}
\newcommand{\beq}{\begin{equation}}
\newcommand{\eeq}{\end{equation}}

\newcommand{\norm}{\mid\mid \cdot \mid\mid}
\newcommand{\nono}{\nonumber}

\newcommand{\lln}{\log\log(n)}
\newcommand{\llnk}{\log\log(n_k)}

\newcommand{\suite}{_{n\geq 1}}
\newcommand{\suik}{_{k\geq 1}}
\newcommand{\Dn}{\Delta_n}
\newcommand{\ovDn}{\overline{\Delta}_n}
\newcommand{\ovDnk}{\overline{\Delta}_{n_k}}
\newcommand{\ovDnkm}{\overline{\Delta}_{n_{k-1}}}
\newcommand{\Dnk}{\Delta_{n_k}}
\newcommand{\DPn}{\Delta\Pi_n}
\newcommand{\DPnk}{\Delta\Pi_{n_k}}

\newcommand{\wtDnk}{\widetilde{\Delta}_{n_k}}

\newcommand{\Dnkm}{\Delta_{n_{k-1}}}
\newcommand{\DFn}{\Delta F_n}
\newcommand{\DPFn}{\Delta\Pi F_n}

\newcommand{\zud}{[0,1]^d}

\begin{document}
\begin{frontmatter}



\title{Some new almost sure results on the functional increments of the uniform empirical process}
\author{Davit Varron}
\ead{dvarron@univ-fcomte.fr}
\address{Université de Franche-Comté, CNRS\thanksref{adresse}}
\thanks[adresse]{Laboratoire de Mathématiques Pures et Appliquées, UMR CNRS 6623, 16, route de Gray,
25030 Besançon Cedex, FRANCE. Phone: +33 3 81 66 65 89}
\begin{abstract}
Given an observation of the uniform empirical process $\alp_n$, its functional increments $\alp_n(u+a_n\cdot)-\alp_n(u)$ can be viewed as a single random process, when $u$ is distributed under the Lebesgue measure. We investigate the almost sure limit behaviour of the multivariate versions of these processes as $\nif$ and $a_n\downarrow 0$. Under mild conditions on $a_n$, a convergence in distribution and functional limit laws are established. The proofs rely on a new extension of usual Poissonisation tools for the local empirical process.
\end{abstract}

\end{frontmatter}

\section{Introduction and statement of the results}
In 1992, Wschebor \cite{Wschebor92} discovered the following property of the Wiener process $\mathbf{W}$ on $[0,1]$ : almost surely, for each $0\le a<b\le 1$, and for each Borel set $B\subset \RRR$, 
\beq \lab\poo \ao u\in [a,b],\;  \e^{-1/2}\po\mathbf{W}(u+\e)-\mathbf{W}(u)\pf \in B\af\pff \mathop{\rar}_{\e\rar 0}(b-a)\PPP\poo \NN(0,1)\in B\pff.\label{pe}\eeq
Here $\lab$ denotes the Lebesgue measure.
That result was later extended to a much wider class of processes by Azaïs and Wschebor \cite{AzaisW96}. It is well known that the increments of the uniform empirical process share several asymptotic behaviors with  the increments of $\mathbf{W}$. Due to the works of many researchers in the past decades, our knowledge on these local functional increments (and also their generalized versions, when indexing by classes of functions, and when the distribution of the $(U_i)_{i\geq 1}$ is non necessarily uniform) is getting deeper and deeper. Strong approximation techniques of these local empirical processes by Gaussian  processes (see \cite{MasonVan}, \cite{EinmahlM4}, \cite{Einmahl2}, \cite{Deheuvels1}) and Poisson processes (see \cite{DeheuvelsM3}, \cite{DeheuvelsM6}, \cite{VarronMaumy1}) have been established, as well as large deviation principles and functional laws of the iterated logarithm of Strassen type (see \cite{DeheuvelsM2}, \cite{Mason2}, \cite{Einmahl3}, \cite{DeheuvelsM5}, \cite{Mason1}) and of so called nonstandard type (see \cite{DeheuvelsM2}, \cite{DeheuvelsM6},\cite{Varron3}, \cite{VarronMaumy1}). Second order results such as clustering rates and Chung-Mogulski laws have been also established (see \cite{Deheuvels2}, \cite{DeheuvelsLif3}, \cite{BerthetLif}, \cite{Berthet1}, \cite{Berthet2}, \cite{Shmileva1}). Of course, the preceding list is non exhaustive, as the complete study of the increments of non functional type (due to the works of Mason, Cs\H{o}rg\H{o}, R\'ev\'esz, Stute among others) was also pioneering in that field. We refer to \cite{MSW}, Chapter 14, and the references therein for an overview on that specific topic. \lb
 A natural question is : can we obtain  results such as (\ref{pe}) for these empirical increments. In this article, we provide a positive answer. In addition, we show that :
\begin{itemize}
\item A similar almost sure convergence in distribution holds for the \textbf{functional} increments of the uniform empirical process (see Theorem \ref{T1});
\item An analogue of the Strassen law of the iterated logarithm for the functional increments of the empirical process (see Mason \cite{Mason2}) holds the present context (see Theorem \ref{loglog});
\item An analogue of the nonstandard functional law for the increments of the empirical process (see \cite{DeheuvelsM3}) also holds (see Theorem \ref{nonstandard});
\item Each of these result hold when handling the increments of the uniform empirical process based on a \textbf{multivariate} sample.
\end{itemize}
Before stating our results, we need to introduce some notations.
Denote by $D(\zud)$ the space of all distribution functions of finite signed measures on  $\zud$, and $\mmi\cdot\mmid$ the sup norm on $\zud$, namely :
$$\mmi f \mmid:=\sup_{t\in \zud} \mid f(t)\mid.$$
For $f\in D(\zud)$ and $A\subset \zud$ Borel, we shall write $f(A)$ for $\mu(A)$, where $\mu$ is the measure associated to $f$.\lb
Consider an i.i.d. sequence $(U_n)\suite$ uniformly distributed on $\zud$.  For $s=(s^{(1)},\ldots,s^{(d)})$ and $t=(t^{(1)},\ldots,t^{(d)})$ belonging to $\RRR^d$ the notation $s\prec t$ shall be understood as $s^{(k)}\le t^{(k)}$ for each $k=1,\ldots,d$. We shall also write the cube $[s,t]:= [s^{(1)},t^{(1)}]\times\ldots \times [s^{(d)},t^{(d)}]$. For fixed $u\in \RRR^d$ and $a\in[0,1]$ we will denote by $u+a$ the vector $(u_1+a,u_2+a,\ldots,u_d+a)$ and define:  
$$\Dn(u,a,\cdot):=\frac{\sliin \poo\mathds{1}_{[u,u+a\cdot ]}(U_i)-\lab([0,a\cdot])\pff}{\sqrt{na^d}}.$$
We shall also write $W$ for the standard Wiener sheet (namely $\Cov(W(t),W(s))$ $:=(s^{(1)}\wedge t^{(1)})\times\ldots\times (s^{(d)}\wedge t^{(d)}$)) and $\lab^*$ (resp. $\lab_*$) the outer (resp. inner) Lebesgue measure on the subsets of $[0,1]^d$. 
Our first result is a multivariate, functional analogue of (\ref{pe}).\begin{theo}\label{T1}
Assume that :
\beq a_n\downarrow 0,\;\;na_n^d\uparrow \infty,\;\;\lin \log(1/a_n)/\log\log(n)>1.\label{conditions cvloi}\eeq
Then almost surely, for each hypercube $I$ fulfilling both $\lab(I)> 0$ and $I\subset [0,1-\dd]^d$ for some $\dd>0$, the following assertions are true :
\begin{align}
\nono (i)& \text{ for each closed set } F\subset D(\zud)\text{ we have } \\
\nono & \lsn \frac{\lab^*\po\{u\in I,\; \Dn(u,a_n,\cdot)\in F\}\pf}{\lab(I)}\le\PPP(W\in F),\\
\nono (ii)& \text{ for each open set } O\subset D(\zud)\text{ we have } \\
\nono & \lin \frac{\lab_*\po\{u\in I,\; \Dn(u,a_n,\cdot)\in O\}\pf}{\lab(I)}\geq \PPP(W\in O).\\
\end{align}
\end{theo}
Our second result is a functional law of the iterated logarithm, in the same vein as Theorem \ref{T1}. We will denote by $J$ the rate function  related to the large deviation properties of a Wiener sheet
\beq  J(f):=\inf\aooo \ili_{\zud} g^2(u) du,\; f=\ili_{[0,\cdot]} g(s) ds\afff,\; f\in D(\zud)\label{J},\eeq
with the convention $\inf_{\emptyset}=+\infty$.
 Definition (\ref{J}) enables us to write the Strassen ball as
\beq \SSS:=\aoo f\in D(\zud),\;J(f)\le 1\aff.\eeq
\begin{theo}\label{loglog}
Assume that :
\beq a_n\downarrow 0,\;\;na_n^d\uparrow \infty,\;\;\frac{na_n^d}{\log\log(n)}\rar \infty,\;\;\lin \frac{\log(1/a_n)}{\log\log(n)}>2.\label{conditions cvloi}\eeq
Then almost surely, for each hypercube $I$ fulfilling both $\lab(I)>0$ and $I\subset [0,1-\dd]^d$ for some $\dd>0$ we have :
\beq \frac{\lab\pooo\aoo u\in I,\; \frac{\Dn(u,a_n,\cdot)}{\sqrt{2\lln}}\leadsto \SSS\aff\pfff}{\lab(I)}=1.\eeq
Here $f_n\leadsto \SSS$ means that the sequence $(x_n)\suite$ has cluster set $\SSS$ in the Banach space $D(\zud)$.
\end{theo}
Our third result is a nonstandard functional law of the iterated logarithm, when $na_n^d\sim c\log\log(n)$ for a constant $0<c<\infty$. To state this, we shall introduce the following rate function ruling the large deviations of a standard Poisson process on $\RRR^d$.
\beq \mathfrak{J}(f):=\inf\aooo \ili_{\zud} h(u) du,\; f=\ili_{[0,\cdot]} g(s) ds\afff,\; f\in D(\zud)\label{mathfrakJ},\eeq
with $h(x):=x\log(x)-x+1$ for $x>0$ and $h(0):=0$. For a constant $c>0$ we shall write
$$\Gam_c:=\aoo f \in D\po \zud\pf,\; \mathfrak{J}(f)\le 1/c\aff.$$
\begin{theo}\label{nonstandard}
Assume that $na_n^d \sim c\log\log(n)$ for some $0<c<\infty$. Then almost surely, for each hypercube $I$ fulfilling both $\lab(I)>0$ and $I\subset [0,1-\dd]^d$ for some $\dd>0$ we have
\beq \frac{\lab\pooo\aoo u\in I,\; \frac{\DFn(u,a_n,\cdot)}{c \log\log(n)}\leadsto \Gam_c\aff\pfff}{\lab(I)}=1,\eeq
where
\beq \DFn(u,a_n,t):=\sliin \mathds{1}_{[u,u+a_n t]}(U_i),\; u,t\in \zud.\label{DFn}\eeq
\end{theo}
In each of our proofs, we systematically use two kinds of key arguments :
\begin{itemize}
\item A tool for replacing probabilities involving the studied processes by probabilities involving their \textit{poissonised} versions. These Poissonised versions have a property that play the same role as the independence of increments (which plays a crucial role in the result of Wschebor \cite{Wschebor92}).
\item The existing knowledge of the asymptotic behavior of probabilities for a \textbf{single} sequence of functional increments (for example the Poissonised version of $\Dn(0,a_n,\cdot)$).
\end{itemize}
For the proof of Theorem 1 (see \S \ref{PT1}), we only use existing results. In particular, we make use of a "Poissonisation" tool of Giné \textit{et. al.}
 (see \S \ref{SP}). For the proofs of Theorems  \ref{loglog} and \ref{nonstandard} (see \S \ref{Ploglog}  and \S \ref{Pnonstandard} respectively), we need an extended version of the just mentioned Poissonsisation tool, which allows us to handle maximal inequalities for sums of i.i.d. processes (those inequalities playing a crucial role in the proofs of the functional laws of Mason \cite{Mason2} and Deheuvels and Mason \cite{DeheuvelsM3}). This extended version is stated and proved in \S \ref{SP}.
 \section{An extended poissonisation tool}\label{SP}
Whenever possible, substituting empirical processes by their Poissonised versions can be very handy, due to the main property of Poisson measures, which can be seen as a generalization of independence of increments for real indexed processes. The following fact, due to Giné-Mason-Zaitsev is, to the best of our knowledge, the most general form of such a Poissonisation lemma, for which the early versions go at least to Einmahl \cite{Einmahl86}.
\begin{fact}[Giné-Mason-Zaitsev, \cite{GineMZ1}, Lemma 2.1]\label{pois}
Let $(D,\DD)$ be a measurable semigroup, $X_0\equiv 0$ and $(X_i)_{i\geq 1}$ be a sequence of $\DD$-measurable, independent, identically distributed random variables. Let $\eta$ be a Poisson variable with expectation $n$ independent of $(X_i)_{i\geq 1}$ and let $B,C\in \DD$ be such that $\PPP(X_1\in B)\le 1/2$. Then
$$\PPP\poo \sliin \mathds{1}_B(X_i)X_i\in C\pff\le 2 \;\PPP\poo \sli_{i=1}^{\eta} \mathds{1}_B(X_i)X_i\in C\pff.$$
As a consequence, for any positive measurable function $H$ we have
$$\EEE\poo H\poo\sliin \mathds{1}_B(X_i)X_i\pff\pff\le 2\;\EEE\poo H\poo \sli_{i=1}^\eta \mathds{1}_{B}(X_i)X_i\pff\pff.$$
\end{fact}
That fact is crucial in our proof of Theorem \ref{T1}. To prove Theorem \ref{loglog}, we shall need an analogue of Fact \ref{pois} for probabilities related to \textbf{maximal inequalities} for partial sums in a Banach space. This analogue is indeed a consequence of a much wider generalization of Fact \ref{pois}, for which we need to introduce some notations.
Given a semigroup $D$ we shall write $\wtD:=\bculi_{n\geq 1}D^n$. Also, given a set $\chi$ we call a \textit{truncating} application any function $\phi: \wtD\mapsto \chi$ for which, for any $p\geq 2$ and $d_1,\ldots,d_{p}\in D^{p}$ we have $\phi(d_1,\ldots,d_p,d_p)= \phi(d_1,\ldots,d_p)$ and $\phi(d_1,d_1,d_2,\ldots,d_p)= \phi(d_1,d_2,\ldots,d_p)$. We shall say that $\phi$ is \textit{zero-irrelevant} when we add the property $\phi(0,d_1,\ldots,d_p)=\phi(d_1,\ldots,d_p)$. We shall write, for simplicity of notations, 
\begin{align}
\nono \sli_{i=q}^{\rar p}d_i:=&\po d_q,d_q+d_{q+1},\ldots,\sli_{i=q}^pd_i\pf, \text{ when }p\geq q,\\
\nono :=& 0\text{ otherwise. }
\end{align}
when $p\geq q$ and $0\in D$ otherwise.
\begin{propo}\label{bonus} Endow $\wtD$ with the $\sig$-algebra $\wt{\DD}:=\ao \wtC\subset \wtD,\;\forall n\geq 1,\;\wtC\cap D^n\in \DD^{\bigotimes n}\af$. Let $(\chi,\AAA)$ be a measurable space and $\phi:\; \po\wtD,\wt{\DD}\pf\mapsto \po\chi, \AAA\pf$ a measurable truncating application. For any $B\in \DD,\;C\in \AAA$ such that $\PPP\po X\in B\pf\le 1/2$ we have
$$\PPP\pooo \phi \poo \sli_{i=1}^{\rar n}\mathds{1}_{B}(X_i)X_i\pff\in C\pfff\le 2 \;\PPP\pooo \phi \poo \sli_{i=1}^{\rar \eta}\mathds{1}_B(X_i)X_i\pff\in C\pfff.$$
As a consequence, for any positive measurable function $H$ we have 
$$\EEE\poooo H\poo \phi \poo \sli_{i=1}^{\rar n}\mathds{1}_{B}(X_i)X_i\pff\pfff\pffff\le 2 \;\EEE\poooo H\pooo \phi \poo \sli_{i=1}^{\rar \eta}\mathds{1}_B(X_i)X_i\pff\pfff\pffff.$$
\end{propo}
\textbf{Proof of Proposition \ref{bonus}: }\lb
To prove that proposition we shall follow the line of the proof of Fact \ref{pois} (see \cite{GineMZ1}, Lemma 2.1) and go a one step further. Write $p_B:=\PPP(X_1\in B)$ and denote by $(\tau_i,Y_i)_{i\geq 1}$ an i.i.d. sequence for which $Y_i$ is independent of $\tau_i$, $\PPP(\tau_i=1)=1-\PPP(\tau_i=0)=p_B$ and $Y_i$ has the distribution of $X_i$ conditionally to $X_i\in B$. A simple calculation shows that $\mathds{1}_B(X_i)X_i=_d \tau_i Y_i$, from where, by conditioning on $(\tau_1,\ldots,\tau_n)$:
\begin{align}
\nono 
\PPP\pooo \phi \poo \sli_{i=1}^{\rar n} \mathds{1}_B(X_i)X_i\pff\in C\pfff=& \sli_{\PP\subset \{1,\ldots,n\}}p_B^{\sharp \PP}(1-p_B)^{n-\sharp \PP}\;\PPP\pooo \phi \poo \sli_{i=1}^{\rar n}\mathds{1}_\PP(i)Y_i\pff\in C\pfff,
\end{align}
with the notation $\sharp \PP$ for the number of elements of $\PP$.\lb
Now notice that, for fixed $\PP$ and for each permutation of indices $\sig$, the law of $(\mathds{1}_\PP(1)Y_1,\ldots,\mathds{1}_\PP(n)Y_n)$ and $(\mathds{1}_\PP(\sig(1))Y_1,\ldots, \mathds{1}_\PP(\sig(n))Y_n)$ are identical. By choosing $\sig$ such that $(\mathds{1}_\PP(\sig(1)),\ldots, \mathds{1}_\PP(\sig(n)))=(1,\ldots,1,0,\ldots,0)$ the vector $\sli_{i=1}^{\rar n} \mathds{1}_{\PP}(\sig(i))Y_i$ has his last $n-\sharp \PP+1$ coordinates equal, from where, since $\phi$ is truncating :
$$\PPP\pooo\phi \poo \sli_{i=1}^{\rar n} \mathds{1}_B(X_i)X_i\pff\in C\pfff=\sli_{k=0}^n {n\choose k} p_B^k(1-p_B)^{n-k}\;\PPP\pooo\phi\poo \sli_{i=1}^{\rar k}Y_i\pff\in C\pfff.$$ 
The remainder of the calculus follows exactly as in the proof of Lemma 2.1 in \cite{GineMZ1}, until the last line, where it suffices to prove that 
\beq\phi\poo \sli_{i=1}^{\rar \eta}\tau_i Y_i\pff=_d\phi \poo \sli_{i=1}^{\rar \eta_B}Y_i\pff,\label{taul}\eeq
where $\eta_B$ is Poisson with expectation $np_B$, independent of $\YYN$.
This is done by writing 
\begin{align}
\nono& \PPP\pooo \phi\poo \sli_{i=1}^{\rar \eta}\tau_i Y_i\pff\in C\pfff\\
\nono =& \sli_{m\geq 0}\PPP(\eta=m)\sli_{\PP\subset\{1,\ldots,m\}}p_B^{\sharp \PP}(1-p_B)^{m-\sharp \PP}\;\PPP\pooo\phi \poo \sli_{i=1}^{\rar m}\mathds{1}_\PP(i)Y_i\pff\in C\pfff\\
\nono=&\sli_{m \geq 0}\frac{n^m}{m!}e^{-n}\sli_{k=0}^m{m\choose k}p_B^k(1-p_B)^{m-k}\;\PPP\pooo \phi \poo \sli_{i=1}^{\rar k} Y_i\pff\in C\pfff\\\nono &\text{ (by the same arguments as above)}\\
\nono=&\sli_{k\geq 0}\PPP\pooo \phi \poo \sli_{i=1}^{\rar k} Y_i\pff\in C\pfff\;\sli_{m\geq k}\frac{n^m}{m!}\frac{m!}{k!(n-k)!}p_B^k(1-p_B)^{m-k}e^{-n}\\
\nono =& \sli_{k\geq 0}\PPP\pooo \phi \poo \sli_{i=1}^{\rar k} Y_i\pff \in C\pff  e^{-n}\frac{(np_B)^k}{k!}\sli_{m'\geq 0}\frac{{(n(1-p_B))}^{m'}}{m'!}\\
\nono =&\sli_{k\geq 0}\PPP(\eta_B=k)\PPP\pooo \phi \poo \sli_{i=1}^{\rar k} Y_i\pff\in C\pfff,
\end{align}
which proves Proposition \ref{bonus}.$\Box$\lb
Our next proposition shows that, if $\phi$ is also zero-irrelevant, then $\phi\poo \sli_{i=1}^{\rar \eta} \mathds{1}_B(X_i)X_i\pff,\;B\in \DD$ has a property of independence which is very similar to the property of independence of Poissonised sums.
\begin{propo}\label{indep}Assume now that $\phi$ is truncating and zero-irrelevant. In the setting of Proposition \ref{bonus}, without imposing that $\EEE(\eta)=n$, if $B_1,B_2,\ldots,B_r$ are disjoint, then
$$\coooo\phi\poo \sli_{i=1}^{\rar \eta} \mathds{1}_{B_1}(X_i)X_i\pff,\ldots, 
\phi\poo \sli_{i=1}^{\rar \eta} \mathds{1}_{B_r}(X_i)X_i\pff\cffff\text{ are mutually independent. }$$
\end{propo}
\textbf{Proof of Proposition \ref{indep}}:\lb
Write $\lab:=\EEE(\eta)$, $p_\ell:=\PPP(X_1\in B_\ell),\;\ell=1,\ldots,r$, $p_{r+1}:=1-\sli_{\ell=1}^r p_\ell$ (assuming without loss of generality that each of these quantities is nonzero) and consider arbitrary events $C_1,\ldots,C_r$.
Now define
\begin{itemize}
\item For each $i\geq 1$, a mutually independent sequence $(\tau_{i1},\ldots,\tau_{ir})_{i\geq 1}$ for which $(\tau_{i1},\ldots,\tau_{ir})=_d (\mathds{1}_{B_1}(X_i),\ldots,\mathds{1}_{B_r}(X_i))$.
\item A mutually independent family $\po Y_{i\ell}\pf_{i\geq 1,\; \ell=1,\;\ldots,r}$, where the $Y_{i\ell}$ are respectively distributed as $X_i\mid X_i\in B_\ell$. 
\item The above-mentioned family are independent from each other.
\end{itemize} 
Direct computations show that, for fixed $i\geq 1$.
$$(\mathds{1}_{B_1}(X_i)X_i,\ldots,\mathds{1}_{B_r}(X_i)X_i)=_d(\tau_{i1}Y_{i1},\ldots,\tau_{ir}Y_{ir}).$$
We have by conditioning successively with respect to $\eta$ and $(\tau_{i,\ell})_{i\geq 1,\ell=1,\ldots,r}$:
\begin{align}
\nono &\PPP\pooo \bcali_{\ell=1}^r\aoo\phi\poo \sli_{i=1}^{\rar \eta} \mathds{1}_{B_\ell}(X_i)X_i\pff\in C_\ell\aff\pfff\\
\nono =& \PPP\poo \bcali_{\ell=1}^r\aoo \phi\poo \sli_{i=1}^{\rar\eta}\tau_{i\ell}Y_{i\ell}\pff\in C_\ell\aff\pff\\
\nono =&\sli_{m\geq 0} \PPP(\eta=m)\PPP\pooo \bcali_{\ell=1}^r\aoo \phi\poo\sli_{i=1}^{\rar m}\tau_{i\ell}Y_{i\ell}\pff\in C_\ell\aff\pfff\\
\nono =&\sli_{m\geq 0} \PPP(\eta=m)\mathop{\sli_{\PP_1\cup\ldots\cup\PP_r\subset \{1,\ldots,m\}}}_{\PP_1,\ldots,\PP_r\text{ disjoint }}\proli_{\ell=1}^{r+1}p_\ell^{\sharp \PP_\ell}\; \PPP\pooo \bcali_{\ell=1}^r\aoo\phi\poo \sli_{i=1}^{\rar m}\mathds{1}_{\PP_\ell}(i)Y_{i\ell}\pff\in C_\ell\aff\pfff,\end{align}
with $\PP_{r+1}:=\{1,\ldots,m\}-\bculi_{\ell=1}^r\PP_\ell$ in the preceding formula.
Let us focus on a single term of the last sum. Writing $k_\ell:=\sharp \PP_\ell,\; \ell=1,\ldots,r+1$, we can find an permutation of indices $\sig$ such that the first $k_1$ integers of $\{1,\ldots,m\}$ are $\sig(i),\;i\in \PP_1$, the next $k_2$ integers are $\sig(i),\;i\in \PP_2$, and so on. As $\phi$ is both truncating and zero-irrelevant, we have almost surely for each $\ell\le r$ (writing $k_0:=0$)
$$\phi \poo \sli_{i=1}^{\rar m} \mathds{1}_{\PP_\ell}(\sig(i))Y_{i\ell}\pff=\phi \poo \sli_{k_0+\ldots+k_{\ell-1}}^{\rar k_0+\ldots+k_\ell} Y_{i\ell}\pff,$$
from where these $r$ random variables are mutually independent. It follows that
\begin{align}
\nono &\PPP\pooo \bcali_{\ell=1}^r\aoo\phi\poo \sli_{i=1}^{\rar \eta} \mathds{1}_{B_\ell}(X_i)X_i\pff\in C_\ell\aff\pfff\\
\nono =& \sli_{m\geq 0}\frac{\lab^m}{m!}e^{-\lab}\mathop{\sli_{k_1,\ldots,k_r\in\{0,\ldots,m\}}}_{k_{r+1}:=m-k_1-\ldots-k_r\geq 0}m!\proli_{\ell=1}^{r+1}\frac{p_\ell^{k_\ell}}{k_\ell!}\proli_{\ell=1}^r\PPP\pooo \phi\poo \sli_{i=1}^{\rar k_\ell}Y_{i,\ell}\pff\in C_\ell\pfff\\
\nono =&\sli_{k_1,\ldots,k_{r+1}\in \NNN}\;\proli_{\ell=1}^{r+1}\frac{(\lab p_\ell)^{k_\ell}}{k_\ell !}e^{-\lab p_\ell}\proli_{\ell=1}^r\PPP\pooo \phi\poo \sli_{i=1}^{\rar k_\ell}Y_{i,\ell}\pff\in C_\ell\pfff.
\end{align}
Now writing $\eta_1,\ldots,\eta_r$ as independent Poisson random variables with respective expectations $\lab p_1,\ldots,\lab p_r$, which are also independent of $Y_{i,\ell},\ell=1,\ldots,r,i\geq 1$, the last expression is equal to 
\begin{align}
\nono\proli_{\ell=1}^r \PPP\pooo \phi\poo \sli_{i=1}^{\rar \eta_\ell}Y_{i\ell}\pff\in C_{\ell}\pfff\sli_{k_{r+1}\geq 0}\frac{(\lab p_{r+1})^{k_{r+1}}}{k_{r+1}!}e^{-\lab p_{r+1}}=\proli_{\ell=1}^r \PPP\pooo \phi\poo \sli_{i=1}^{\rar \eta_\ell}Y_{i\ell}\pff\in C_{\ell}\pfff.
\end{align} 
The proof is concluded by applying (\ref{taul}) with the formal replacement of $\eta_B,Y_{i},\tau_{i}$ by $\eta_\ell,Y_{i\ell},\tau_{i\ell}$ for $\ell=1,\ldots,r$.$\Box$\lb
\textbf{Remarks}: First, notice that Fact \ref{pois} can be deduced from Proposition \ref{bonus} with the choice of $\phi(d_1,\ldots,d_p):=d_p$.\lb
Second, note that, given a collection of functions $\rho_{\ell}: D\rar \RRR$ (with $\ell=1,\ldots,r$) the application 
$$\phi:(d_1,\ldots,d_p)\rar \coo\max_{i=1,\ldots,p}\;\rho_1(d_i),\ldots,\max_{i=1,\ldots,p}\;\rho_r(d_i)\cff $$ is truncating. Moreover, if each $\rho_\ell$ attains its minimum at $0\in D$, the application $\phi$ is zero-irrelevant. For particular choices of $\rho$, we readily obtain two results that may have an interest in themselves.\lb
The choice of $\rho:=\pm\mathds{1}_C$ leads to the following corollary.
\begin{coro}
Under the setting of Fact \ref{pois} we have : 
$$\PPP\poo \exists m\le n,\; \sliim \mathds{1}_{B}(X_i)X_i\in C\pff\le 2\; \PPP\poo \exists m\le \eta,\; \sliim \mathds{1}_{B}(X_i)X_i\in C\pff,$$
$$\PPP\poo \forall m\le n,\; \sliim \mathds{1}_{B}(X_i)X_i\in C\pff\le 2\; \PPP\poo \forall m\le \eta,\; \sliim \mathds{1}_{B}(X_i)X_i\in C\pff.$$
Now, dropping the assumption $\EEE(\eta=n)$ and taking  $B_1,\ldots,B_r,C_1,\ldots, C_r\in \DD$ with $B_1,\ldots B_r$ disjoint we have:
\begin{enumerate}
\item If $0\notin \bculi_{\ell=1}^rC_\ell$, then 
$$\PPP\poo \bcali_{\ell=1}^r\aoo\exists m\le \eta,\; \sliim \mathds{1}_{B_{\ell}}(X_i)X_i\in C_\ell\aff\pff=\proli_{\ell=1}^r\PPP\poo \exists m\le \eta,\; \sliim \mathds{1}_{B_{\ell}}(X_i)X_i\in C_\ell\pff,$$
\item If $0\in \bcali_{\ell=1}^rC_\ell$, then $$\PPP\poo \bcali_{\ell=1}^r\aoo\forall m\le \eta,\; \sliim \mathds{1}_{B_{\ell}}(X_i)X_i\in C_\ell\aff\pff=\proli_{\ell=1}^r\PPP\poo \forall m\le \eta,\; \sliim \mathds{1}_{B_{\ell}}(X_i)X_i\in C_\ell\pff.$$
\end{enumerate}
\end{coro}
Next, the choice of $\rho(\cdot)$ as a  semi norm leads to 
\begin{coro}\label{tep}
Under the setting of Fact \ref{pois} we have, if $(D,\norm)$ is a semi normed space for which $\norm$ is $\DD$ measurable :
$$\PPP\poo \max_{m\le n}\Mmi \sli_{i=1}^m \mathds{1}_B(X_i)X_i\Mmi \in C\pff\le 2\;\PPP\poo \max_{m\le \eta}\Mmi \sli_{i=1}^m \mathds{1}_B(X_i)X_i\Mmi \in C\pff,$$
$$\EEE\pooo H\poo \max_{m\le n}\Mmi \sli_{i=1}^m \mathds{1}_B(X_i)X_i\Mmi\pff\pff\le2\; \EEE\pooo H\poo \max_{m\le \eta}\Mmi \sli_{i=1}^m \mathds{1}_B(X_i)X_i\Mmi\pff\pfff.$$
Now dropping the assumption $\EEE(\eta=n)$ and taking with $B_1,\ldots B_r$ disjoint, the random variables
$$\coo\max_{m\le\eta}\Mmi \sliim \mathds{1}_{B_1}(X_i)X_i\Mmi,\ldots,\max_{m\le\eta}\Mmi \sliim \mathds{1}_{B_r}(X_i)X_i\Mmi\cff$$
are mutually independent.
\end{coro}
Note that the last statement of Corollary \ref{tep} can be deduced more straightforwardly by making use of the independence properties of Poisson point processes.\lb
Roughly speaking, the preceding corollary shows that blocking arguments for partial sums of i.i.d. random variables can be Poissonised (for which we still have independence properties of Poisson measures). In our proof of Theorem \ref{loglog} we shall use the particular function $\phi$ defined as follows. The semigroup $D$ will be taken to be $D([0,1]^d)^{[0,1]^d}$, $\chi:=[0,\infty)^{\zud}$ and
$$\phi(d_1,\ldots,d_p):=\coo\max_{i=1,\ldots,p}\;\mmi d_i(u)\mmid\cff_{u\in \zud}.$$
\section{Proof  of Theorem \ref{T1}}\label{PT1}
Choose $\dd>0$ and a hypercube $I\subset[0,1-\dd]^d$ for which $\lab(I)>0$ . We can assume without loss of generality that $\lab(I)<1/2$. By a finite union argument, the full version of Theorem \ref{T1} shall readily follow.\lb
Making use of the usual tools in the theory of weak convergence in $D(\zud)$ (see, e.g., \cite{Vander}, Chapter 1.5) together with standard arguments of countable union/intersection of events, we need to establish the following proposition (in what follows we write $\mid u\mid_d$ for $\max\{\mid u_k\mid,\;k=1,\ldots,d\}$).
\begin{propo}\label{fo}
For each integer $p\geq 1$, $\tht_1,\ldots,\tht_p\in \RRR^p$ and $t_1,\ldots,t_p\in \zud$ we have
\beq \ili_I\exp\poo i\sli_{j=1}^p\tht_j\Dn(u,a_n, t_j)\pff du\cvps \lab(I)\exp\po -\frac{1}{2}\tht \Sig \tht\pf,\label{fidi}\eeq
where $\Sig[k,k']:=\Cov(W(t_k),W(t_{k'}))=\lab([0,t_k]\cap [0,t_{k'}])$.\lb
For each $\e>0$, there exists $\dd>0$ such that, almost surely
\beq \lsn\lab\pooo\aoo u\in I,\;\sup_{\mid s-t\mid_d<\dd}\mid \Dn(u,a_n, s)-\Dn(u,a_n, t)\mid>\e\aff\pfff\le \e.\label{oscillations}\eeq
\end{propo}
To prove this proposition, we shall apply the following result. For measurability concerns, we shall endow the space $D(\zud)$ with the $\sig$-algebra $\TT$ spawned by the applications : $$P_{t_1,\ldots,t_p}(f):=(f(t_1),\ldots,f(t_p)),\;p\geq 1,\;t_1,\ldots,t_p\in \zud.$$
Clearly $\TT$ coincides with the $\sig$ algebra spawned by the balls related to the norm $\norm_{\zud}$. We shall also consider the Poissonised version of $\Dn(\cdot,\cdot)$, namely
\beq \DPn(u,a,\cdot):=\frac{\sli_{i=1}^{\eta_n} \po \mathds{1}_{[u,u+a\cdot]}(U_i) -\lab([0,a\cdot])\pf}{\sqrt{na^d}},\label{DPn}\eeq
where $\eta_n$ is a Poisson random variable with expectation $n$ and independent of $\UUN$.  Our proof of Proposition \ref{fo} relies on the following Proposition.
\begin{propo}\label{variance}Let $\rho_n$ be a sequence of measurable applications from $(D(\zud),\TT)$ to $\CCC$. Then
\begin{align}
\nono&\EEE\poooo \pooo\frac{1}{\lab(I)}\ili_I\rho_n(\Dn(u,a_n,\cdot))du-\EEE\poo \rho_n\po\DPn(0,a_n,\cdot)\pf\pff\pfff^2\pffff\\
\nono=&O(a_n)\Var\poo \rho_n(\DPn(0,a_n,\cdot))\pff.
\end{align}
\end{propo} 
\textbf{Proof:}
We shall apply Fact \ref{pois}. We choose the semigroup $D$ to be the space $D(\zud)$, endowed with the $\sig-$algebra $\DD:=\TT$.
Clearly the applications of the form
$$\Psi: f\rar \ili_I\rho\po f([u,u+a\cdot])\pf du,\;a\in[0,1],\;\rho \text{ measurable from }(D(\zud),\TT)\text{ to }\CCC.$$
are $\TT$ measurable.\lb
We take $B:=\ao f\in D(\zud),f(I+[0,a_n]^d)+\lab\po I+[0,a_n]^d\pf>0\af$ and $X_i:=
\mathds{1}_{[0,\cdot]}(U_i)-\cdot$. Clearly, $X_i$ are all $\TT$ measurable and $\PPP(X_1\in B)=\PPP(U_1\in I+[0,a_n]^d)\le 1/2$ (for all large $n$). We then consider the applications
\begin{align}
\nono &H_n:f\rar \pooo\ili_I\frac{1}{\lab(I)}\rho_n\poo \frac{f([u,u+a_n\cdot])}{\sqrt{na_n^d}}\pff du-\EEE\poo \rho_n\po\DPn(0,a_n,\cdot)\pf\pff\pfff^2,\\
\nono &\text{ which satisfies, for all }\mathfrak{n}\geq 1\\
\nono &H_n\poo \sli_{i=1}^\mathfrak{n} \mathds{1}_B(X_i)X_i\pff=_{a.s.}\;H_n\poo \sli_{i=1}^\mathfrak{n} X_i\pff.
\end{align}
Applying Fact \ref{pois} for fixed $n$ leads to the bound
\begin{align}
\nono&\EEE\poooo \pooo\frac{1}{\lab(I)}\ili_I\rho_n(\Dn(u,a_n,\cdot))du-\EEE\poo \rho_n\po\DPn(0,a_n,\cdot)\pf\pff\pfff^2\pffff\\
\nono \le & 2\EEE\poooo \pooo\frac{1}{\lab(I)}\ili_I\rho_n(\DPn(u,a_n,\cdot))du-\EEE\poo \rho_n\po\DPn(0,a_n,\cdot)\pf\pff\pfff^2\pffff\\
\nono =&2\Var \pooo\frac{1}{\lab(I)}\ili_I\rho_n(\DPn(u,a_n,\cdot))du\pfff\text{, as soon as }I+[0,a_n]^d\subset \zud,\\
\nono =&2\ili_I\ili_I\Cov\pooo\frac{1}{\lab(I)}\ili_I\rho_n(\DPn(u,a_n,\cdot)),\frac{1}{\lab(I)}\ili_I\rho_n(\DPn(v,a_n,\cdot))\pfff du dv
\end{align} 
Now, for fixed $u,v$ satisfying $[u,u+a_n]\cap[v,v+a_n]=\emptyset$ , the corresponding covariance is null, as $\DPn(u,a_n,\cdot)\;\indep \;\DPn(v,a_n,\cdot)$ (this can be seen for example by choosing $B_1:=\ao f\in D(\zud),\; f([u,u+a_n])+\lab([u,u+a_n])>0\af$, $B_2:=\ao f\in D(\zud),\;f([v,v+a_n])+\lab([v,v+a_n])>0\af$ and $\phi(d_1,\ldots,d_p):=d_p$ and applying Proposition \ref{indep}). This entails :
\begin{align}
\nono &\ili_I\ili_I\Cov\pooo\frac{1}{\lab(I)}\ili_I\rho_n(\DPn(u,a_n,\cdot)),\frac{1}{\lab(I)}\ili_I\rho_n(\DPn(v,a_n,\cdot))\pfff du dv\\
\nono =& \mathop{\ili_{u,v\in I^2}}_{\mid u-v\mid_d\le a_n}\Cov\pooo\frac{1}{\lab(I)}\ili_I\rho_n(\DPn(u,a_n,\cdot)),\frac{1}{\lab(I)}\ili_I\rho_n(\DPn(v,a_n,\cdot))\pfff du dv \\
\nono \le & \frac{1}{\lab(I)^2}\mathop{\ili_{u,v\in I^2}}_{\mid u-v\mid_d\le a_n}\sqrt{\Var\poo\ili_I\rho_n(\DPn(u,a_n,\cdot))\pff}\sqrt{\Var\poo\ili_I\rho_n(\DPn(v,a_n,\cdot))\pff} du dv\\
\nono =&O\pooo a_n\Var\poo\ili_I\rho_n(\DPn(0,a_n,\cdot))\pff\pfff. \Box\end{align}
\textbf{Proof of Proposition \ref{fo}}
The following fact shall be needed to prove Proposition \ref{fo}. To the best of our knowledge, it has not yet be written in the literature. However, it can be readily proved by making use of modern tools in empirical processes theory.
\begin{fact}\label{cvpoisson}
The sequence $\DPn(0,a_n,\cdot)$ converges in distribution to $W$ in the space $\po D(\zud),\norm\pf$.
\end{fact}
The proof of Proposition \ref{fo} will be achieved in two steps.\lb
\textbf{Step 1: proof of Proposition \ref{fo} along a subsequence}\lb
Consider the subsequence 
\beq n_k:=\coo \exp\po k/\log(k)\pf\cff,\label{nk}\eeq so that  $\sli_{k\geq 1}\ank<\infty$ by assumption (\ref{conditions cvloi}). Here $[a]$ stands for the unique integer $m$ fulfilling $m\le a<m+1$. By taking, for arbitrary $p\geq 1,\;t_1,\ldots,t_p\in \zud$ and $\tht_1,\ldots,\tht_p\in \RRR$ the function
$$\rho_n:=f\rar \exp\poo i\sli_{j=1}^p\tht_j f(t_j)\pff,$$
we readily obtain by Proposition \ref{variance} that, almost surely
$$\limk \ili_I \exp\poo i\sli_{j=1}^p\tht_j\Dnk(u,\ank,t_j)\pff du-\lab(I)\EEE\poo \exp\poo i\sli_{j=1}^p\tht_j\DPnk(0,\ank,t_j)\pff\pff=0,$$
from where point (\ref{fidi}) of Proposition \ref{fo} is proved along $(n_k)\suik$, making use of Fact \ref{cvpoisson}.\lb
We now fix $\e>0$ and choose, for fixed $\dd>0$ :
 $$\rho_n:\;f\in D\po \zud\pf \rar \mathds{1}_{(\e,\infty)}\poo\sup_{\mid t'-t\mid_d<\dd}\mid f(t')-f(t)\mid\pff.$$
 Then, by Proposition \ref{variance} :
 \begin{align}
 \nono \EEE\poooo\pooo&\frac{1}{\lab(I)}\lab\poo\aoo u \in I,\;\sup_{\mid t'-t\mid_d<\dd}\mid \Dn(u,a_n,t')-\Dn(u,a_n, t)\mid >\e\aff\pff\\
 &\;-\PPP\poo \sup_{\mid t'-t\mid_d<\dd}\mid\DPn(0,a_n,t')-\DPn(0,a_n,t)\mid\pff\pfff^2\pffff=O(a_n).
 \end{align}
Again, Proposition \ref{cvpoisson} entails, almost surely :
\begin{align}
\nono &\limk \lab\poo\aoo u \in I,\;\sup_{\mid t'-t\mid_d<\dd}\mid \Dnk(u,\ank,t')-\Dnk(u,\ank, t)\mid >\e\aff\pff\\
=&\PPP\poo\sup_{\mid t'-t\mid_d<\dd}\mid W(t')-W(t)\mid>\e\pff,\label{te1}
\end{align}
which proves point (\ref{oscillations}) of Proposition \ref{fo} along $(n_k)\suik$, as $W$ admits a uniformly continuous version on $\zud$.\lb
\textbf{Step 2: Blocking arguments}\lb
Now consider the block $N_k:=\{\nkm+1,\ldots,\nk\}$. As $\nk/\nkm\rar 1$ and $\ank/\ankm\rar 1$, we just need to prove (\ref{fidi}) by replacing $\Dn(\cdot,a_n,\cdot)$ by
\beq \ovDn(u,s):=\frac{\sliin \poo \mathds{1}_{[u,u+\ank v]}(U_i)-\lab([0,\ank s])\pff}{\sqrt{n_k\ank^d}},\;k\geq 1,\;n\in N_k,\label{ovDn}\eeq
which satisfies $\ovDnk(\cdot,\cdot)=\Dnk(\cdot,\ank,\cdot)$ almost surely for each $k\geq 1$. Notice that, for $n\in N_k,\;p\geq 1,\;t_1,\ldots,t_p\in \zud,\;\tht_1,\ldots,\tht_p\in \RRR$ :
\begin{align}& \Mid \ili_I\exp\poo i\sli_{j=1}^p \tht_j\ovDn(u,t_j)\pff du-\ili_I\exp\poo i\sli_{j=1}^p \tht_j\ovDnk(u,t_j)\pff du\Mid\\
\nono \le & \max_{j=1,\ldots,p}\mid \tht_j\mid \pooo\ili_I \Mmi \ovDn(u,\cdot)-\ovDnk(u, \cdot)\Mmid du\pfff.
\end{align}
Moreover, for fixed $\e>0$ and $\delta>0$ we have almost surely :
\begin{align}
\nono &\lab\poo \aoo u\;:\sup_{\mid t'-t\mid_d\le \dd}\Mid \ovDn(u,t')-\ovDn(u,t)\Mid>4\e\aff\pff\\
\nono \le & \lab\poo \aoo u\;:\sup_{\mid t'-t\mid_d\le \dd}\Mid \ovDnk(u,t')-\ovDnk(u,t)\Mid>\e\aff\pff\\
\nono&+\lab\poo \aoo u\;:\Mmi \ovDnk(u,\cdot)-\ovDn(u,\cdot)\Mmid>\e\aff\pff,
\end{align}
where the almost sure limit of the first term is known by (\ref{te1}). It turns out that the proof of Proposition \ref{fo} shall be completed if we can show that
\beq \max_{n\in N_k}\ili_I \Mmi \ovDn(u,\cdot)-\ovDnk(u, \cdot)\Mmid du\cvps 0.\label{te2}\eeq
By making use of the Montgomery-Smith maximal inequality (see \cite{Montgom}, Theorem 1 and Corollary 4), we know that, for fixed $\e>0$ :
\begin{align}
\nono&\PPP\poo \max_{n\in N_k}\ili_I \Mmi \ovDn(u,\cdot)-\ovDnk(u, \cdot)\Mmid du>30\e\pff\\
\nono \le &9\;\PPP\pooo \ili_I \Mmi \ovDnk(u,\cdot)-\ovDnkm(u, \cdot)\Mmid du>\e\pff\\
 =&9\;\PPP\poo \ili_I\Mmi \sli_{i=1}^{\nk-\nkm}\mathds{1}_{[u,u+\ank \cdot]}(U_i)-\lab([0,\ank \cdot])\Mmid du> \e\sqrt{\nk\ank^d}\pff.\label{bound}
\end{align}
These probabilities shall be controlled as follows.
\begin{lem}\label{u1}
We have, $\eta_{\nk-\nkm}$ denoting Poisson random variable with expectation $\nk-\nkm$, independent of $(U_i)_{i\geq 1}$ :
\begin{align}
&\EEE\pooo \ili_I\Mmi \sqrt{\frac{\nk-\nkm}{\nk}}\Delta\Pi_{\nk-\nkm}(u,\ank ,\cdot)\Mmi du\pfff\rar 0\label{p1},\\
&\nono\EEE\poooo\frac{\nk-\nkm}{\nk}\pooo \ili_I\Mmi \Delta_{\nk-\nkm}(u,\ank, \cdot)\Mmi du -\EEE\poo\ili_I\Mmi \Delta\Pi_{\nk-\nkm}(u,\ank, \cdot)\Mmi du\pff\pfff^2\pffff\\
=&O(\ank).\label{p2}
\end{align}
\end{lem}
\textbf{Proof:} The second point is a straightforward adaptation of Proposition \ref{variance}, while the first point comes from the fact that, for all large $k$:
\begin{align}
\nono&\EEE\pooo \ili_I\Mmi \sqrt{\frac{\nk-\nkm}{\nk}}\Delta\Pi_{\nk-\nkm}(u,\ank, \cdot)\Mmid du\pfff\\
\nono =&\lab(I)\times\EEE\pooo\Mmi \sqrt{\frac{\nk-\nkm}{\nk}}\Delta\Pi_{\nk-\nkm}(0,\ank, \cdot)\Mmid\pfff\\
=& \lab(I)\times\sli_{m\geq 0}\frac{(\nk-\nkm)^m}{m!}\exp\po -(\nk-\nkm)\pf \mu_m,\label{ddd}
\end{align}
where, $\mu_0:=0$ and, for each $m\geq 1$:
\begin{align}
\nono \mu_m:=&\EEE\poo\Mmi \sqrt{\frac{m}{\nk}}\Delta_{m}(0,\ank, \cdot)\Mmid\pff\\
 \le &C_0\sqrt{\frac{m}{\nk}}\label{dd},
\end{align}
where $C_0$ is a universal constant.
Note that (\ref{dd}) can be proved by a bracketing numbers argument. For example, apply Corollary 19.35, p. 288 in \cite{VanderAsymptotic} with $\FF:=\ao \mathds{1}_{[0,a_nt]},\;t\in \zud\af$, $F=\mathds{1}_{[0,a_n]^d}$ and $P$ the uniform distribution on $\zud$.\lb
Inserting the bound (\ref{dd}) in (\ref{ddd}) yields
\begin{align}
\nono\EEE\pooo \ili_I\Mmi \sqrt{\frac{\nk-\nkm}{\nk}}\Delta\Pi_{\nk-\nkm}(u,\ank, \cdot)\Mmid du\pfff \le & \frac{C_0}{\sqrt{\nk}}\EEE\poo\sqrt{\eta_{\nk-\nkm}}\pff\\
\nono=&O\pooo\sqrt{\frac{\nk-\nkm}{\nk}}\pfff\\
\nono =&o(1).\Box
\end{align}
Now we can prove (\ref{te2}) by taking an arbitrary $\e>0$, applying the bound (\ref{bound}), then making use of point (\ref{p1}) of Lemma \ref{u1} to obtain, for all large $k$ : 
\begin{align}
\nono &\PPP\pooo \ili_I\Mmi \sli_{i=1}^{\nk-\nkm}\mathds{1}_{[u,u+\ank \cdot]}(U_i)-\lab([0,\ank \cdot])\Mmi du> \e\sqrt{\nk\ank^d}\pfff\\
\nono=&\PPP\pooo \ili_I\Mmi \sqrt{\frac{\nk-\nkm}{\nk}}\Delta_{\nk-\nkm}(u,\ank, \cdot)\Mmi du >\e\pfff\\
\nono \le& \PPP \pooo \Mid \ili_I\Mmi \sqrt{\frac{\nk-\nkm}{\nk}}\Delta_{\nk-\nkm}(u,\ank, \cdot)\Mmi\\
\nono &\;\;-\EEE\poo \ili_I\Mmi \sqrt{\frac{\nk-\nkm}{\nk}}\Delta\Pi_{\nk-\nkm}(u,\ank, \cdot)du\Mmi\pff\Mid>\e/2\pfff.
\end{align}
We then apply point (\ref{p2}) of Lemma \ref{u1} together with Markov's Inequality and the Borel-Cantelli Lemma.$\Box$
\section{Proof of Theorem \ref{loglog}}\label{Ploglog}
We first need a large deviation result for $\DPn(0,a_n,\cdot)$. We shall write (recalling (\ref{J}))
\begin{align}
\nono J(A):=&\inf\aoo J(f),\;f\in A\aff,\;A\subset D(\zud).
\end{align}
\begin{propo}\label{gdpoisson}
Under the assumptions $a_n\rar 0$ and $na_n/\log\log(n)\rar \infty$ we have :
\begin{itemize}
\item For each closed set $F\in \TT$ of $\po D(\zud),\norm\pf$:
$$\lsn \frac{1}{\log\log(n)}\log \poooo\PPP\pooo \frac{\DPn(0,a_n,\cdot)}{\sqrt{2\log\log(n)}}\in F\pfff\pffff\le -J(F),$$
$$\lsn \frac{1}{\log\log(n)}\log \poooo\PPP\pooo \frac{\Dn(0,a_n,\cdot)}{\sqrt{2\log\log(n)}}\in F\pfff\pffff\le -J(F).$$
\item For each open set $O\in \TT$ of $\po D(\zud),\norm\pf$:
$$\lin \frac{1}{\log\log(n)}\log \poooo\PPP\pooo \frac{\DPn(0,a_n,\cdot)}{\sqrt{2\log\log(n)}}\in O\pfff\pffff\geq -J(O),$$
$$\lin \frac{1}{\log\log(n)}\log \poooo\PPP\pooo \frac{\Dn(0,a_n,\cdot)}{\sqrt{2\log\log(n)}}\in O\pfff\pffff\geq -J(O).$$
\end{itemize}
\end{propo}
\textbf{Proof :} The part concerning $\Dn(0,a_n,\cdot)$ is a consequence of Proposition 3.2 in \cite{VarronISUP}. The proof of the part concerning $\DPn(0,a_n,\cdot)$ is very similar to the proof of Proposition 1 in \cite{Mason1}. We omit details. $\Box$\lb
We can assume without loss of generality that $\lab(I)<1/2$. The proof shall be split in two parts.  
\subsection{Upper bounds}
This subsection is devoted to proving that, almost surely :
\beq \lab\poo\bcali_{\e>0\in \QQQ}\;\bcali_{n_0\geq 1}\;\bculi_{n\geq n_0}\aoo u\in I,\; \frac{\Dn(u,a_n,\cdot)}{\sqrt{2\lln}}\notin \SSS^\e\aff\pff=0,\label{upper}\eeq
where $\SSS^\e=\aoo f\in D(\zud),\; \inf\ao \mmi f-g\mmid,\;g\in\SSS \af<\e\aff.$\lb
\textbf{Step 1: proof along a subsequence}\lb
Take $(\nk)\suik$ as in (\ref{nk}). For fixed $\e>0$ we shall show that there exists $\dd>0$ for which, almost surely as $\kif$ :
\beq \lab\poo\aoo u\in I,\frac{\Dnk(u,\ank,\cdot)}{\sqrt{2\llnk}}\notin \SSS^\e\aff\pff=O\poo\exp(-(1+\dd)\llnk)\label{bu}\pff.\eeq
This achieve this, fist notice that, as $\SSS$ is compact and $J$ is lower semi continuous on $\po D(\zud),\norm_{[0,1]^d}\pf$ we can choose $\dd>0$ so as $J\poo D(\zud)-\SSS^\e\pff>1+2\dd$. Moreover, as $$\lin\log(1/a_n)/\log\log(n)>2,$$ we can assume without loss of generality that $a_n\log(n)^{2+2\dd}\rar 0$. We then make use  of Proposition \ref{variance} with :
$$\phi_n(f):=\exp\po(1+\dd)\lln\pf\mathds{1}_{D(\zud)-\SSS^\e}\poo \frac{f}{\sqrt{2na_n^d\lln}}\pff,$$
which yields
\begin{align}
\nono &\EEE\poooo \pooo e^{(1+\dd)\lln}\ili_{u\in I} \mathds{1}_{D(\zud)-\SSS^\e}\poo \frac{\Dn(u,a_n,\cdot)}{\sqrt{2\lln}}\pff\\
\nono &\;\;\;\;\;\;-e^{(1+\dd)\lln}\PPP\poo \frac{\DPn(0,a_n,\cdot)}{\sqrt{2\lln}}\notin \SSS^\e\pff\pfff^2\pffff\\
\nono =&O(a_n)\PPP\poo \frac{\DPn(0,a_n,\cdot)}{\sqrt{2\lln}}\notin \SSS^\e\pff e^{2(1+\dd)\lln}.
\end{align}
From Proposition \ref{gdpoisson} we know that the last quantity is $o(a_n)\exp\po (1+\dd)\lln\pf=o\po\log(n)^{-1-\dd}\pf$ which is sumable along $\nk$. This proves (\ref{bu}) and also proves (\ref{upper}) along $(\nk)\suik$.\lb
\textbf{Step 2 : blocking arguments}\lb
Now take $\ovDn(\cdot,\cdot)$ as defined in (\ref{ovDn}).
We shall now show that, almost surely 
\beq \frac{\lab\aooo u \in I,\; \limk \max_{n\in N_k}\frac{\mmi \ovDn(u,\cdot)-\ovDnk(u,\cdot)\mmi}{\sqrt{2\llnk}}=0\afff}{\lab(I)}=1.\label{oscps}\eeq
For fixed $k\geq 1$ we shall apply Proposition \ref{bonus} in the following setting: we take semigroup $D:=D([0,1]^d)^{[0,1]^d}$,  endowed with $\DD:= \TT^{\bigotimes [0,1]}$. We take $\chi:=[0,\infty)^{\zud}$ and
$$\phi(d_1,\ldots,d_p):=\coo\max_{i=1,\ldots,p}\;\mmi d_i([u,u+\ank \cdot])\mmid\cff_{u\in \zud}.$$ We apply Proposition \ref{bonus} to the sequence $X_m:=\mathds{1}_{[0,\cdot]}(U_m)-\cdot,\; m\geq \nkm+1$, with $n:=\nk-\nkm$ and  
\begin{align}
\nono H:\;& g\;\rar \;\pooo \frac{\ili_{u\in I} \mathds{1}_{A_k}\po g(u)\pf du}{\lab(I)}-\mathfrak{m}_k\pff^2\text{, where}\\
\nono A_k:=& [\sqrt{2\nk\ank^d\llnk}\e,+\infty),\\
\nono \mathfrak{m}_k:=&\EEE\pooo \frac{1}{\lab(I)}\ili_I \mathds{1}_{A_k}\poo\max_{m\le \eta_{\nk-\nkm}} \Mmi \sli_{i=\nkm+1}^{\nkm+m} \mathds{1}_{[u,u+\ank,\cdot]}(U_i)-\lab([0,\ank\cdot])\Mmid\pff du \pfff.
\end{align}
Clearly, if we take the event $B_k:=\aoo f\in D(\zud),\; f\po I+[0,\ank]^d\pf+\lab\po I+[0,\ank]^d\pf>0\aff$ we have $\PPP(X_i\in B)=\PPP\po U_i\in I+[0,\ank]^d\pf<1/2$ (for $k$ large enough), and, for all $\mathfrak{n}\geq 1$:
$$H\poo\phi\poo\sli_{\nkm+1}^{\rar \nkm+\mathfrak{n}}\mathds{1}_{B_k}(X_i)X_i\pff\pff=_{a.s.}\;H\poo\phi\poo\sli_{\nkm+1}^{\rar \nkm+\mathfrak{n}}X_i\pff\pff,$$
From where, by Proposition \ref{bonus} (writing $v_k:= \e\sqrt{2\nk\ank^d\llnk}$)
\begin{align}
\nono &\EEE\poooo \pooo \frac{\lab\poo \aoo u\in I,\; \mathop{\max}_{n\in N_k}\mmi \ovDn(u,\cdot)-\ovDnk(u,\cdot)\mmid>\e\sqrt{2\llnk}\aff\pff}{\lab(I)}-\mathfrak{m}_k\pfff^2\pfff\\
\nono \le & \frac{2}{\lab(I)^2}\Var\pooo \ili_I \mathds{1}_{[v_k,+\infty)}\poo\max_{m\le \eta_{\nk-\nkm}} \Mmi \sli_{i=\nkm+1}^{\nkm+m} \mathds{1}_{[u,u+\ank,\cdot]}(U_i)-\lab([0,\ank\cdot])\Mmid\pff du \pfff\\
\nono =& \frac{2}{\lab(I)^2}\ili_{I^2}\Cov\pooo\mathds{1}_{[v_k,+\infty)}\poo\max_{m\le \eta_{\nk-\nkm}} \Mmi \sli_{i=\nkm+1}^{\nkm+m} \mathds{1}_{[u,u+\ank,\cdot]}(U_i)-\lab([0,\ank\cdot])\Mmid\pff,\\
\nono &\;\mathds{1}_{[v_k,+\infty)}\poo\max_{m\le \eta_{\nk-\nkm}} \Mmi \sli_{i=\nkm+1}^{\nkm+m} \mathds{1}_{[v,v+\ank,\cdot]}(U_i)-\lab([0,\ank\cdot])\Mmid\pff\pfff du dv.
\end{align}
For fixed $v$ fulfilling $\mid s-t\mid>\ank$, the corresponding covariance is null, as the two involved random variables are independent. To see this, apply Proposition \ref{indep} with $B_1:=\ao f\in D(\zud),\; f([u,u+\ank])+\lab([u,u+\ank])>0\af$ and $B_2:=\ao f\in D(\zud),\; f([v,v+\ank])+\lab([v,v+\ank])>0\af$. As $\sli_k\ank<\infty$, assertion (\ref{oscps}) will be proved as soon as we prove that 
\beq\sli_{k\geq 1}\PPP \poo\max_{m\le \eta_{\nk-\nkm}} \Mmi \sli_{i=\nkm+1}^{\nkm+m} \mathds{1}_{[0,\ank,\cdot]}(U_i)-\lab([0,\ank\cdot])\Mmid\geq v_k\pff <\infty.\label{sumable}\eeq
To show (\ref{sumable}) we split the probabilities in two
\begin{align}
\nono \PPP_k:=&\PPP\poo\max_{m\le \eta_{\nk-\nkm}} \Mmi \sli_{i=\nkm+1}^{\nkm+m} \mathds{1}_{[u,u+\ank,\cdot]}(U_i)-\lab([0,\ank\cdot])\Mmid\geq v_k\pff\\
\nono \le & \PPP\poo \max_{m\le 2(\nk-\nkm)} \Mmi \sli_{i=\nkm+1}^{\nkm+m} \mathds{1}_{[u,u+\ank,\cdot]}(U_i)-\lab([0,\ank\cdot])\Mmid\geq v_k\pff\\
\nono &+\PPP\poo\eta_{\nk-\nkm}>2(\nk-\nkm)\pff.
\end{align}
The second term is sumable by the Chernoff bound. The first term can be bounded by a sequence of the form 
$$s_k:=\exp(-A\frac{\nk}{\nk-\nkm}\llnk)+\exp(-B\sqrt{\nk\ank^d\llnk}),$$ by making use an inequality of Talagrand (see Inequality A.1 in \cite{EinmahlM1}), with $M=1$, $\GG:=\ao \mathds{1}_{[0,\ank s]},\;s\in \zud\af$ and $t= \frac{\e}{2}\sqrt{2\nk\ank^d\llnk}$, together with the following first moment bound for symetrised empirical processes:
\begin{align}
\nono &\EEE\pooo \Mmi \e_i\sli_{i=\nkm+1}^{\nkm+2(\nk-\nkm)}\mathds{1}_{[u,u+\ank,\cdot]}(U_i)\Mmid\pfff\\
\nono\le &2\EEE\pooo \Mmi \sli_{i=\nkm+1}^{\nkm+2(\nk-\nkm)}\mathds{1}_{[u,u+\ank,\cdot]}(U_i)-\lab([0,\ank\cdot])\Mmid\pfff\\
\le &2C_0\sqrt{(\nk-\nkm)\ank^d\llnk}\pff\label{r}\\
\nono =& o(\sqrt{\nk\ank^d\llnk}),
\end{align}
where (\ref{r}) is no more than the already proved inequality (\ref{dd}). As $\llnk\sim\log(k),\; (\nk-\nkm)/\nk\rar 0,\;\nk\ank^d/\llnk\rar \infty$, we deduce that $\sli_k s_k<\infty$, which implies that $\sli_k\PPP_k<\infty$.
\subsection{Lower bounds}
By compactness of $\SSS$, we just need to prove that, for fixed $f\in \SSS,\;n_0\in \NNN$ and $\e>0$, almost surely:
\beq \lab\poo \bcali_{n\geq n_0}\;\ao u\in I,\; \Dn(u,a_n,\cdot)\notin f^\e\af\pff=0.\label{low}\eeq
\textbf{Step 1: proof for a modified sequence}\lb
As $f\in \SSS$, we have $J(f^\e)=1-2\dd$ for some $\dd>0$. Consider the subsequence $n_k:=\coo \exp\po k^{1+\dd}\pf\cff$ and write
$$\wtDnk(u,t):=\frac{\sli_{i=\nkm+1}^{\nk}\mathds{1}_{[u,u+\ank t]}(U_i)-\lab([0,\ank t])}{\sqrt{2(\nk-\nkm)\ank^d \llnk}},$$
which defines a sequence of mutually independent processes.
As $(\nk-\nkm)\ank^d\sim \nk\ank^d$, we can use  Proposition (\ref{gdpoisson}) and obtain, for all large $k$:
\beq\PPP\poo \wtDnk(0,\cdot)\in f^\e\pff \geq  \exp\po -(1-\dd)\llnk\pf\geq\frac{1}{k^{1-\dd^2}},\eeq
from where, when $m\rar \infty$ :
$$m^{\dd^2}=O\poooo\sli_{k=1}^m\PPP\poo \wtDnk(0,\cdot)\in f^\e\pff\pffff.$$
As the preceding  events are mutually independent we obtain for a constant $C$: 
\begin{align}
\nono \PPP\poo \bcali_{k=1}^m \aoo \wtDnk(0,\cdot)\notin f^\e\aff\pff\le& \exp\pooo -\PPP\poo \wtDnk(0,\cdot)\in f^\e\pff\pfff\\
\nono \le  \exp \po - Cm^{\dd^2}\pf.
\end{align}
Hence by Markov's inequality we have, for fixed $\tau>0$  and $k_0$ large enough to fulfill $I+[0,a_{n_{k_0}}]^d\subset \zud$, and $m\geq k_0$:
\begin{align}
\nono &\PPP\pooo \lab\poo \bcali_{k=k_0}^m\ao u\in I,\; \wtDnk(u,\cdot)\notin f^\e\af\pff>\tau\pfff\\
\nono &\PPP\poo \ili_{u\in I}\proli_{k=k_0}^m\mathds{1}_{D(\zud)-f^\e}\po \wtDnk(u,\cdot)\pf>\tau\pff\\
\nono \le & \frac{\lab(I)}{\tau}\PPP\poo \bcali_{k=k_0}^m \aoo \wtDnk(0,\cdot)\in f^\e\aff\pff\\
\nono \le & \frac{\lab(I)}{\tau}\exp \po -Cm^{\dd^2}\pf,
\end{align}
which is sumable in $m$, from where we obtain that, almost surely as $m\rar \infty$ :
$$\lab\poo \bcali_{k=k_0}^m\aoo u\in I,\;  \wtDnk(u,\cdot)\notin f^\e\aff\pff \rar 0,$$
whence, with probability one :
\beq \lab\poo \bculi_{k\geq 1}\bcali_{k\geq k_0}\aoo u\in I,\; \wtDnk(u,\cdot)\notin f^\e\aff\pff=0.\label{low}\eeq
\textbf{Step 2 : proof for the original sequence}\lb
In view of (\ref{low}), and since
\begin{align}
\nono\frac{\nk-\nkm}{\nk}\rar &1,\\
\nono \frac{\Dnk(\cdot,\ank,\cdot)}{\sqrt{2\llnk}}=&\sqrt{\frac{\nkm}{\nk}}\frac{\Dnkm(\cdot,\ank,\cdot)}{\sqrt{2 \llnk}}+\sqrt{\frac{\nk-\nkm}{\nk}}\wtDnk(\cdot,\cdot),
\end{align}
we just need to show that, almost surely, as $k_0\rar \infty$ :
\beq \lab\poo \bculi_{k\geq k_0}\aoo\sqrt{\frac{\nkm}{\nk}}\Mmi\frac{\Dnkm(u,\ank,\cdot)}{\sqrt{2\llnk}}\Mmid>\e\aff\pff\rar 0.\label{lu}\eeq
Again, by Markov's inequality we get, for fixed $\tau>0$ and $k_0$ large enough :
\begin{align} 
\nono &\PPP\pooo \lab\poo \bculi_{k\geq k_0}\aoo\sqrt{\frac{\nkm}{\nk}}\Mmi\frac{\Dnkm(u,\ank,\cdot)}{\sqrt{2 \llnk}}\Mmid>\e\aff\pff\geq \tau\pfff\\
\nono \le &\frac{\lab(I)}{\tau}\sli_{k=k_0}^\infty \PPP\pooo\sqrt{\frac{\nkm}{\nk}}\Mmi\frac{\Dnkm(0,\ank,\cdot)}{\sqrt{2\llnk}}\Mmid>\e\pfff.
\end{align}
Now applying Markov's inequality once again, together with the bound (\ref{dd}) we get that
\begin{align}
\nono &\sli_{k=k_0}^\infty \PPP\pooo\sqrt{\frac{\nkm}{\nk}}\Mmi\frac{\Dnkm(0,\ank,\cdot)}{\sqrt{2\llnk}}\Mmid>\e\pfff\\
\nono \le & \frac{1}{\e}\sli_{k=k_0}^\infty\EEE\pooo\sqrt{\frac{\nkm}{\nk}}\Mmi\frac{\Dnkm(0,\ank,\cdot)}{\sqrt{2\llnk}}\Mmid\pfff\\
\nono =&O\poo\sli_{k=k_0}^\infty\sqrt{\frac{\nkm}{\nk\llnk}}\pff,
\end{align}
which is sumable in $k_0$. This proves (\ref{lu}) and hence completes the proof of Theorem \ref{loglog}.
\section{Proof of Theorem \ref{nonstandard}}\label{Pnonstandard}
The proof is in the same vein as the proof of Theorem \ref{loglog}. Hence, to avoid lengthy redundancies we shall only focus on the single technical difference (even if the real novelty of the present proof relies on the non-written methods that mimic the proof of Theorem \ref{loglog}). First, we shall require the following result, which is included in Proposition 3.2. in \cite{VarronMaumy1}. Here $\DPFn(\cdot,a_n,\cdot)$ stands for the poissonised version of $\DFn(\cdot,a_n,\cdot)$, namely:
$$\DPFn(u,a_n,t):=\sli_{i=1}^{\eta_n}\mathds{1}_{[u,u+a_n t]}(U_i).$$
\begin{fact}\label{gdp}
Under the assumption $na_n^d\sim c\log\log(n)$ we have
\begin{itemize}
\item For each closed set $F\in \TT$ we have
$$\lsn \frac{1}{c\log\log(n)}\log\poo\PPP \poo \frac{\DPFn(0,a_n,\cdot)}{c\log\log(n)} \in F\pff\pfff \le -\mathfrak{J}(F);$$ 
\item For each open set $O \in \TT$ we have
$$\lin \frac{1}{c\log\log(n)}\log\pooo\PPP \poo \frac{\DPFn(0,a_n,\cdot)}{c\log\log(n)} \in O\pff \pfff\geq -\mathfrak{J}(O).$$ 
\end{itemize}
\end{fact}
The proof of Theorem \ref{nonstandard} is achieved following the same steps as in the proof of Theorem \ref{loglog}, replacing $\DPn(\cdot,a_n,\cdot)/\sqrt{2n\lln}$ by $\DPFn(\cdot,a_n,\cdot)/c\lln$, and $\SSS$ by $\Gam_c$. The only point where the methodology changes is when proving an analogue of (\ref{sumable}), namely $\sli \PPP'_k<\infty$, where
$$\PPP'_k:=\PPP\poo \max_{m \le \eta_{\nk-\nkm}}\Mmi \sli_{i=\nkm+1}^{\nkm+m}\mathds{1}_{[0,\ank\cdot]}(U_i)\Mmi \geq c\e \llnk\pff.$$
 But, in this particular case, as all the sumed processes on $\zud$ are distribution functions of positive measures, we have
 $$\PPP'_k=\PPP\poo \sli_{i=1}^{\eta_{\nk-\nkm}}\mathds{1}_{[0,\ank]^d}(U_i)\geq c\e \llnk\pff,$$
 where the involved random variable is Poisson with expectation $(\nk-\nkm)\ank^d=o(\llnk)$. Hence we avoid making use of Talagrand's inequality (which does not provide a strong enough bound when $na_n^d\sim c\log\log(n)$), and just apply the Chernoff bound for Poisson random variables to establish that $\sli \PPP'_k<\infty.$ $\Box$ \lb
 \textbf{Aknowledgements}: The author would like to thank the referees for their very useful comments.

\begin{thebibliography}{10}

\bibitem{AzaisW96}
J.M. Azais and M.~Wschebor.
\newblock {Almost Sure Oscillation of Certain Random Processes}.
\newblock {\em Bernoulli}, 2(3):257--270, 1996.

\bibitem{Berthet1}
P.~Berthet.
\newblock {On the rate of clustering to the Strassen set for increments of the
  empirical process}.
\newblock {\em J. Theoret. Probab.}, 10(3):557--579, 1997.

\bibitem{Berthet2}
P.~Berthet.
\newblock {Inner rates of coverage of Strassen type sets by increments of the
  uniform empirical and quantile processes}.
\newblock {\em Stochastic Process. Appl.}, 115(3):493--537, 2005.

\bibitem{BerthetLif}
P.~Berthet and M.~Lifshits.
\newblock Some exact rates in the functional law of the iterated logarithm.
\newblock {\em Ann. Inst. H. Poincaré Probab. Statist.}, 38:811--824, 2002.

\bibitem{Deheuvels2}
P.~Deheuvels.
\newblock Chung type functional laws of the iterated logarithm for tail
  empirical processes.
\newblock {\em Ann. Inst. H. Poincaré Probab. Statist.}, 36:583--616, 2000.

\bibitem{Deheuvels1}
P.~Deheuvels.
\newblock {Strong approximation of quantile process by iterated Kiefer
  processes}.
\newblock {\em Ann. Probab.}, 28(2):909--945, 2000.

\bibitem{DeheuvelsLif3}
P.~Deheuvels and M.~Lifshits.
\newblock {Probabilities of hiting of shifted small balls by centered Poisson
  processes}.
\newblock {\em J. Math. Sci. (N. Y.)}, 118(6):5541--5554, 2003.

\bibitem{DeheuvelsM3}
P.~Deheuvels and D.M. Mason.
\newblock Nonstandard functional laws of the iterated logarithm for tail
  empirical and quantile processes.
\newblock {\em Ann. Probab.}, 18:1693--1722, 1990.

\bibitem{DeheuvelsM2}
P.~Deheuvels and D.M. Mason.
\newblock Functional laws of the iterated logarithm for the increments of
  empirical and quantile processes.
\newblock {\em Ann. Probab.}, 20:1248--1287, 1992.

\bibitem{DeheuvelsM5}
P.~Deheuvels and D.M. Mason.
\newblock Functional laws of the iterated logarithm for local empirical
  processes indexed by sets.
\newblock {\em Ann. Probab.}, 22:1619--1661, 1994.

\bibitem{DeheuvelsM6}
P.~Deheuvels and D.M. Mason.
\newblock Nonstandard local empirical processes indexed by sets.
\newblock {\em J. Statist. Plann. Inference}, 45:91--112, 1995.

\bibitem{Einmahl86}
J.H.J. Einmahl.
\newblock Multivariate empirical processes.
\newblock Technical report, {Ph.D dissertation, Katolieke Universiteit
  Nigmegen}, {1986}.

\bibitem{Einmahl3}
J.H.J. Einmahl.
\newblock The as behavior of the weighted empirical process and the lil for the
  weighted tail empirical process.
\newblock {\em Ann. Probab.}, 20(2):681 -- 695, 1992.

\bibitem{Einmahl2}
J.H.J. Einmahl.
\newblock Extensions of results of koml\'os, major, and tusn\'ady to the
  multivariate case.
\newblock {\em Stochastic Process. and their Appl.}, 79(1):31 -- 58, 1997.

\bibitem{EinmahlM4}
U.~Einmahl and D.M. Mason.
\newblock Gaussian approximation of local empirical processes indexed by
  functions.
\newblock {\em Probab. Theory Related Fields}, 107(3):283--311, 1997.

\bibitem{EinmahlM1}
U.~Einmahl and D.M. Mason.
\newblock Poisson and gaussian approximation of weighted local empirical
  processes.
\newblock {\em J. Theoret. Probab.}, 13:1--13, 2000.

\bibitem{GineMZ1}
E~Giné, D.M. Mason, and A.~Zaitsev.
\newblock {The $ L_1 $ -norm density estimator process}.
\newblock {\em Ann. Probab.}, 31:719--768, 2003.

\bibitem{Mason2}
D.M. Mason.
\newblock A strong invariance principle for the tail empirical process.
\newblock {\em Ann. Inst. H. Poincaré Probab. Statist.}, 24:491--506, 1988.

\bibitem{Mason1}
D.M. Mason.
\newblock A uniform functional law of the iterated logarithm for the local
  empirical process.
\newblock {\em Ann. Probab.}, 32(2):1391--1418, 2004.

\bibitem{MasonVan}
D.M. Mason and W.R. van Zwet.
\newblock {A refinement of the KMT inequality for the uniform empirical
  process}.
\newblock {\em Ann. Probab.}, 15:871--884, 1987.

\bibitem{VarronMaumy1}
M.~Maumy and D.~Varron.
\newblock Non standard functional limit laws for the increments of the compound
  empirical distribution function.
\newblock {\em To appear in Electron. J. Stat.}, 

\bibitem{Montgom}
J.S. Montgommery-Smith.
\newblock Comparison of sums of identically distributed random vectors.
\newblock {\em {Probab. Math. Statist.}}, 14:281--285, 1993.

\bibitem{Shmileva1}
E.~Yu. Shmileva.
\newblock Small ball probabilities for a centered poisson process of high
  intensity. (russian).
\newblock {\em  Zap. Nauchn. Sem. S.-Peterburg. Otdel. Mat. Inst. Steklov.
  (POMI)}, 298:280--303, 2003.

\bibitem{MSW}
G.R. Shorack and J.A. Wellner.
\newblock {\em {Empirical Processes and applications to statistics}}.
\newblock Springer, 1986.

\bibitem{VanderAsymptotic}
A.W. Van~der Vaart.
\newblock {\em Asymptotic Statistics}.
\newblock Cambridge University Press, 1998.

\bibitem{Vander}
A.W. Van~der Vaart and J.A. Wellner.
\newblock {\em Weak convergence and empirical processes}.
\newblock Springer, 1996.

\bibitem{VarronISUP}
D.~Varron.
\newblock Some uniform in bandwidth functional results for the tail uniform
  empirical and quantile processes.
\newblock {\em Ann. I.S.U.P.}, 50(1-2):83--103, 2006.

\bibitem{Varron3}
D.~Varron.
\newblock A nonstandard uniform functional limit law for the increments of the
  multivariate empirical distribution function.
\newblock {\em Adv. Appl. Stat. Sci.}, 1(2):399--428, 2010.

\bibitem{Wschebor92}
M.~Wschebor.
\newblock {Sur les accroissements du processus de Wiener}.
\newblock {\em {C.R. Acad. Sci. Paris, Ser. I}}, 315(12):1293--1296, 1992.

\end{thebibliography}

\end{document}